# SEARCHING FOR A TRAIL OF EVIDENCE IN A MAZE

By Ery Arias-Castro,[1] Emmanuel J. Candès,[2]
Hannes Helgason and Ofer Zeitouni[3]

*University of California, San Diego, California Institute of Technology,
California Institute of Technology and University of Minnesota,
Minneapolis and Weizmann Institute, Israel*

Consider a graph with a set of vertices and oriented edges connecting pairs of vertices. Each vertex is associated with a random variable and these are assumed to be independent. In this setting, suppose we wish to solve the following hypothesis testing problem: under the null, the random variables have common distribution $N(0,1)$ while under the alternative, there is an unknown path along which random variables have distribution $N(\mu,1)$, $\mu > 0$, and distribution $N(0,1)$ away from it. For which values of the mean shift $\mu$ can one reliably detect and for which values is this impossible?

Consider, for example, the usual regular lattice with vertices of the form

$$\{(i,j): 0 \leq i, -i \leq j \leq i \text{ and } j \text{ has the parity of } i\}$$

and oriented edges $(i,j) \to (i+1, j+s)$, where $s = \pm 1$. We show that for paths of length $m$ starting at the origin, the hypotheses become distinguishable (in a minimax sense) if $\mu_m \gg 1/\sqrt{\log m}$, while they are not if $\mu_m \ll 1/\log m$. We derive equivalent results in a Bayesian setting where one assumes that all paths are equally likely; there, the asymptotic threshold is $\mu_m \approx m^{-1/4}$.

We obtain corresponding results for trees (where the threshold is of order 1 and independent of the size of the tree), for distributions other than the Gaussian and for other graphs. The concept of the predictability profile, first introduced by Benjamini, Pemantle and Peres, plays a crucial role in our analysis.

Received May 2007; revised May 2007.
[1]Supported in part by NSF Grant DMS-06-03890.
[2]Supported in part by NSF Grants ITR ACI-0204932 and CCF515362.
[3]Supported in part by NSF Grant DMS-05-03775.
*AMS 2000 subject classifications.* Primary 62C20, 62G10; secondary 82B20.
*Key words and phrases.* Detecting a chain of nodes in a network, minimax detection, Bayesian detection, predictability profile of a stochastic process, martingales, exponential families of random variables.

**1. Introduction.** This paper discusses the model problem of detecting whether or not there is a chain of connected nodes in a given network which exhibit an "unusual behavior." Suppose we are given a graph $G$ with vertex set $V$ and a random variable $X_v$ attached to each node $v \in V$. In that sense, this is a graph-indexed process. We observe a realization of this process and wish to know whether all the variables at the nodes have the same behavior in the sense that they are all sampled from a common distribution $F_0$, or whether there is a path in the network, that is, a chain of consecutive nodes connected by edges, along which the variables at the nodes have a different distribution $F_1$. In other words, can one tell whether hidden in the background noise, there is a chain of nodes that stand out?

Suppose, for example, that $F_0$ is the standard normal distribution, whereas $F_1$ is a normal distribution with mean 0.1 and variance 1. In a situation where the number of nodes along the path we wish to detect is comparably small, the largest values of $X_v$ are typically off this path. Can we reliably detect the existence of such a path? More generally, how subtle an effect can we detect? In this paper, we attempt to provide quantitative answers to such questions by investigating asymptotic detection thresholds—values of the mean shift at which detection is possible and values at which detection by any method whatsoever is impossible.

Detection thresholds depend, of course, on the type of graphs under consideration and we propose the study of two representative graphs which are, in some sense, far from each other, as well as emblematic—regular lattices and trees. We introduce them next. Later in the paper, we will also consider other graphs.

- *Regular lattice in dimension* 2. Our first graph is a regular lattice with nodes

  $$V_m = \{(i,j) : 0 \le i \le m-1, -i \le j \le i \text{ and } j \text{ has the parity of } i\}$$

  and with oriented edges $(i,j) \to (i+1, j+s)$, where $s = \pm 1$. We call $(0,0)$ the *origin* of the graph. Here and below, we use the subscript $m$ in $V_m$ to remind the reader of the radius of the graph. A path in the graph is represented in Figure 1.

- *Complete binary tree.* Our second model is the oriented regular binary tree. The nodes in the tree are of the form

  $$V_m = \{(i,j) : 0 \le i \le m-1, 0 \le j < 2^i\}.$$

  and it has oriented edges $(i,j) \to (i+1, 2j+s)$, where $s \in \{0,1\}$. Again, we call $(0,0)$ the origin of the graph and the subscript $m$ indicates the radius of the graph (i.e., the depth of tree). A path in the tree is represented in Figure 2.

Note that even though the numbers of paths of length $m$ in both graphs are the same, the numbers of nodes are considerably different—about $m^2/2$ for the lattice and $2^m$ for the binary tree.

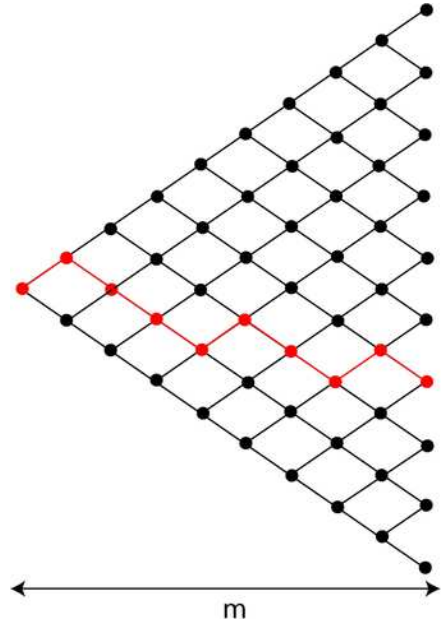

Fig. 1. *Representation of a path (in red) in the regular lattice.*

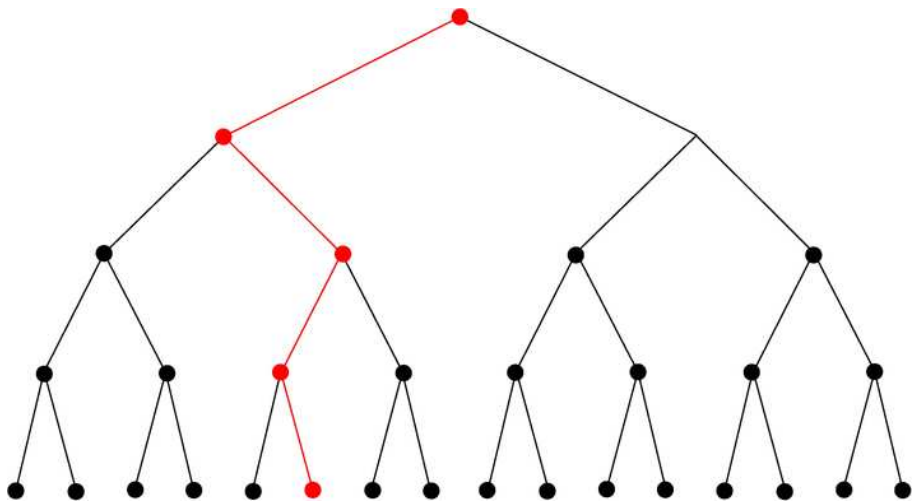

Fig. 2. *Representation of a path (in red) in the binary tree.*

We denote by $\mathcal{P}_m$ the set of paths in the graph starting at the origin and of length $m$. (In this paper, we define the length of a path to be the number of vertices the path visits.) We attach a random variable $X_v$ to each node $v$ in the graph. We observe $\{X_v : v \in V\}$ and consider the following hypothesis testing problem:

- Under $H_0$, all the $X_v$'s are i.i.d. $N(0,1)$.
- Under $H_{1,m}$, all the $X_v$'s are independent; there is an unknown path $p \in \mathcal{P}_m$ along which the $X_v$'s are i.i.d. $N(\mu_m, 1)$, $\mu_m > 0$, while they are i.i.d. $N(0,1)$ away from the path.

In plain English, we would like to know whether there is a path along which the mean is elevated.

1.1. *Motivation.* While this paper is mainly concerned with the study of fundamental detection limits, our problem is in fact motivated by applications in various fields, especially in the area of signal detection.

Suppose we are given very noisy data of the form

$$(1.1) \qquad y_i = S_i + z_i, \qquad i = 1, \ldots, n,$$

where $(S_i)$ are sampled values of a signal of interest and $(z_i)$ is a noise term. Based on the observations $(y_i)$, one would like to decide whether or not a signal is hiding in the noise. That is, we would like to test whether $S = 0$ or not. Suppose, further, that the signal is completely unknown and does not depend on a small number of parameters. In image processing, the signal $S$ might be the indicator function of a general shape we wish to detect or a curve embedded in a two-dimensional pixel array [3]. In signal processing, the signal may be a chirp, a high-frequency wave with unknown and rapidly changing oscillatory patterns [10].

In these situations, we cannot hope to generate a family of candidate signals that would provide large correlations with the unknown signal as the number of such candidates would be exponentially large in the signal size. In response to this obstacle, recent papers [10, 13] have proposed a very different approach, in which the family of candidate signals actually corresponds to a path in a network. We briefly explain the main idea. In most situations, it is certainly possible to generate a family of templates $(\phi_v)_{v \in V}$ which provide good *local* correlations with the signal of interest, for example, over shorter time intervals. Any signal of interest could then be closely approximated by a chain of such templates. Here, a chain is a path in a graph $G$ with nodes $v \in V$ indexed by our templates and rules for connecting templates, these rules possessing the following property: any consecutive sequence of templates in the graph must correspond to a meaningful signal; that is, a signal one might expect to observe (e.g., imagine connecting linear segments to approximate smooth curves). Now, calculate a $Z$-score for each template and denote it $X_v$. For simplicity, assume that $X_v \sim N(\mu_1, 1)$ if the template matches the signal $S$ locally and $X_v \sim N(\mu_0, 1)$ otherwise. Assume $\mu_1 > \mu_0$. Then, the signal detection problem is this: is there a path along which the mean of the $Z$-scores is slightly elevated?

To make things a little more concrete, suppose the unknown signal $S(t)$ is a chirp of the general form $A(t) \exp(i\lambda\varphi(t))$, where $A(t)$ is a smooth amplitude, $\varphi(t)$ is a smooth phase function and $\lambda$ is a large base frequency. Roughly speaking, a chirp is an oscillatory signal with "instantaneous frequency" given by the derivative of the phase, that is, $\lambda\varphi'(t)$. Here, one might use as templates chirplets of the form $\phi_v(t) \propto 1_{I_v}(t) \exp(i(a_v t^2/2 + b_v t))$ which are supported on the time interval $I_v$ and assume the linear instantaneous frequency $a_v t + b_v$. Such templates provide a local quadratic approximation of the unknown phase function $\lambda\varphi(t)$ (or a local linear approximation of the unknown instantaneous frequency) and can exhibit high correlations with the unknown signal, provided that the discretization of the chirplet parameters is sufficiently fine. The chirplet graph [10] then connects pairs of chirplets supported on contiguous time intervals by imposing a certain kind of continuity of the instantaneous frequency in such a way that a path represents a chirping signal with a piecewise linear instantaneous frequency which obeys a prescribed regularity criterion. Given the data vector $y$ (1.1), one would then compute all the chirplet coefficients $X_v = \langle y, \phi_v \rangle$ of $y$. Testing whether there is signal or not amounts to testing whether all the node variables $X_v$ in the chirplet graph have mean 0 or whether there is a path along which the mean is nonzero (the constraint that all possible paths start at a given vertex corresponds to the constraint that if a signal exists, its instantaneous frequency at time 0 is known).

Although the signal detection problem motivates the theoretical study presented in this paper, the problem of detecting a path in a network seems to represent a fundamental abstraction as many modern statistical detection problems can reasonably be formulated in this way. Indeed, it is very easy to imagine that one has available a number of measurements about variables related through a graphical model and that one wishes to detect whether there is a sequence of connected nodes which exhibit a peculiar behavior. We give one example to stimulate the reader's imagination. In [22], water quality in a network of streams is assessed by performing a chemical analysis at various locations along the streams. As a result, some locations are marked as problematic. We may view the set of all tested locations as nodes and connect pairs of adjacent nodes located on the same stream, thereby creating a tree (although not a regular tree), with the root corresponding to the point which is the most downstream. We then assign to each node the value 1 or 0, according to whether the location is problematic or not. A possible model would assume that the variables are Bernoulli, taking the value 1 with probability equal to $p_0$ when the location is normal and $p_1$ when it is anomalous. One can then imagine that one would like to detect a path (or a family of paths) upstream of a certain sensitive location, in order to trace the existence of a polluter, or look for the existence of an anomalous path upstream from the root of the system; see [22]. Note that here, one could also be interested in detecting whether or not there is a family of anomalous paths, as opposed to just one such path. Examples of this kind truly abound; for example, one could imagine detecting atypical gene behaviors in a given gene network, and so on.

1.2. *A quick look at the results.* The optimal detection threshold discussed above is the minimum value of $\mu = \mu_m$ which allows us to reliably tell whether or not there is a path which does not follow the null distribution. This value depends on the criterion used for judging the quality of the decision rule, and statistical decision theory essentially offers two paradigms: the Bayesian and the minimax approach. We study them both.

Consider the minimax paradigm first. Recall that a test $T_m$ is a $\{0,1\}$-valued, measurable function of the collection $(X_v)_{v \in V}$. The minimax risk of a test $T_m$ is defined as

$$(1.2) \quad \gamma(T_m) = \mathbf{P}(\text{Type I}) + \sup_{p \in \mathcal{P}_m} \mathbf{P}_{1,p}(\text{Type II}).$$

Throughout, we write $\mathbf{P}_0$ for the law of $(X_v)$ under $H_0$ and $\mathbf{P}_{1,p}$ for the law of the same variables under $H_{1,m}$ with path $p \in \mathcal{P}_m$. With this notation, Type I and II are shorthand for errors of Type I and II. In longhand,

$$\mathbf{P}(\text{Type I}) = \mathbf{P}_0(T_m = 1), \quad \mathbf{P}_{1,p}(\text{Type II}) = \mathbf{P}_{1,p}(T_m = 0).$$

We say that a sequence of tests $(T_m)$ is *asymptotically powerful* if

$$\lim_{m \to \infty} \gamma(T_m) = 0$$

and *asymptotically powerless* if

$$\liminf_{m \to \infty} \gamma(T_m) \geq 1.$$

When there exists an asymptotically powerful sequence of tests, we say that reliable detection is possible; when all sequences of tests are asymptotically powerless, we say that detection is (essentially) impossible.

1.2.1. *The regular lattice.* We first consider the regular lattice in dimension 2.

THEOREM 1.1. *Consider the regular lattice in dimension 2. Suppose that $\mu_m (\log m)^{1/2} \to \infty$ as $m \to \infty$. That then is a sequence of tests which is asymptotically powerful. On the other hand, suppose that $\mu_m \log m (\log \log m)^{1/2} \to 0$ as $m \to 0$. Every sequence of tests is then asymptotically powerless.*

Theorem 1.1 states that one can detect a path as long as $\mu_m \gg (\log m)^{-1/2}$, while this is impossible if $\mu_m < (\log m)^{-(1+\epsilon)}$ for each $\epsilon > 0$, provided that $m$ is sufficiently large. The reader will note the discrepancy between the lower and the upper bound, which we will comment on in the concluding section.

It turns out that the detection level is radically different in a Bayesian framework where one assumes that all paths are equally likely. For a prior $\pi$ on $\mathcal{P}_m$, namely on paths of length $m$, the corresponding risk of a test $T_m$ is now defined as

$$\gamma_\pi(T_m) = \mathbf{P}(\text{Type I}) + \mathbf{E}_\pi \mathbf{P}_{1,p}(\text{Type II}), \tag{1.3}$$

where $\mathbf{E}_\pi$ stands for the expectation over the prior path distribution, namely, when the path $p$ is drawn according to $\pi$. We adopt the same terminology as before and say that $(T_m)$ is asymptotically powerful if $\gamma_\pi(T_m) \to 0$ and powerless if $\liminf \gamma_\pi(T_m) \geq 1$. The Bayes test associated with $\pi$ is, of course, optimal here. The following theorem shows that under the uniform prior on paths, the optimal Bayesian detectability threshold is about $m^{-1/4}$.

THEOREM 1.2. *Consider the regular lattice in dimension 2 and assume the uniform prior on paths. If $\mu_m m^{1/4} \to \infty$ as $m \to \infty$, then the Bayes test is asymptotically powerful. Conversely, if $\mu_m m^{1/4} \to 0$ as $m \to 0$, then the Bayes test is asymptotically powerless.*

Roughly speaking, if the anomalous path is chosen uniformly at random, one can asymptotically detect it as long as the intensity along the path exceeds $m^{-1/4}$, while no method whatsoever can detect below this level.

Both results indicate that it is possible to detect an anomalous path event when $\mu_m \to 0$ (sufficiently slowly). Note that while one can certainly reliably detect in such circumstances, it may be impossible to tell which sequence of nodes the anomalous path is traversing. This is an example of a situation where detection is possible, but estimation may not be.

1.2.2. *The binary tree.* We are now interested in the complete binary tree.

THEOREM 1.3. *If $\mu_m = \mu \geq \sqrt{2\log 2}$, then there is a sequence of tests that is asymptotically powerful. On the other hand, if $\mu_m = \mu < \sqrt{2\log 2}$, then there is no sequence of tests that is asymptotically powerful. Moreover, if $\mu_m \to 0$ as $m \to \infty$, then every sequence of tests is asymptotically powerless.*

Notice that there is no sharp threshold phenomenon here, in the sense that the minimax risk does not converge to 1 if $\mu_m = \mu < \sqrt{2\log 2}$. For example, the risk of the test which rejects the null hypothesis for large values of the variable at the root node is bounded away from 1 for any $\mu > 0$.

For any graph, and under the normal model, consider the generalized likelihood ratio test (GLRT) which is the test rejecting the null for large values of $M_m := \max\{X_p : p \in \mathcal{P}_m\}$, where $X_p$ is the sum of the node variables along the path $p$:

$$X_p = \sum_{v \in p} X_v. \tag{1.4}$$

The proof of Theorem 1.3 then shows that for the binary tree, the GLRT achieves the minimax threshold in that it has asymptotically full power when $\mu > \sqrt{2 \log 2}$. In this sense, the GLRT rivals the Bayes test under the uniform prior on paths, which, by symmetry, is minimax.

1.3. *Innovations and related work.* In the regular graph model, the number of variables needed to describe the path is $m$, while the total number of nodes or observations is about $m^2/2$. Hence, the topic of this paper fits into the broad framework of nonparametric detection as the object we wish to detect is simply too complex to be reduced to a small number of parameters. Because the theory and practice of detection have been centered around parametric models in which the generalized likelihood ratio test has played a crucial role (see the literature on scan statistics, matched filters and deformable templates, to name a few equivalent terms used in various fields of science and engineering [2, 16, 20, 24]), methods and results for nonparametric detection are comparably underdeveloped. Against this background, we will first provide some evidence showing that the generalized likelihood ratio test does not perform very well in our nonparametric set-up. Our work also differs from the important literature on nonparametric detection in that it does not assume that the unknown object we wish to detect lies in a traditional smoothness class, such as Sobolev or Besov classes, or belongs to an $\ell_p$-ball or some related geometric body; see the book by Ingster and Suslina [19] and the multiple references therein. In fact, our model, techniques and results have nothing to do with this literature and hence our paper contributes to developing the important area of nonparametric detection in what appears to be a new direction. In fact, we are not familiar with statistical theory posing a problem as a graph detection problem and giving precise quantitative bounds. It has come to our attention, however, that Berger and Peres have very recently considered problems which are mathematically closely related to our framework but with a different motivation.

Our paper also has some connections with the theory and practice of multiple hypothesis testing. Indeed, we are interested in situations where testing at each node separately offers little or no power so that we need to combine information from different nodes. Because the anomalous nodes are located on a path, the search naturally involves testing over paths. There are many such paths, however, and in this sense, our problem resembles that of testing many hypotheses (one hypothesis test would be whether the mean along a specified path is zero or not).

1.4. *Organization of the paper.* The paper is organized as follows. In Section 2, we study the detection problems over the regular lattice in dimension 2 and prove our results about the minimax and Bayesian detection thresholds, namely, Theorems 1.1 and 1.2. In Section 3, we prove the detection thresholds for the binary tree. In Section 4, we extend our results to exponential distributions at the nodes and in Section 5, to other distributions and other graphs. In Section 6, we report on numerical simulations which complement our theoretical study. Finally, we conclude with Section 7, where we comment on our findings and discuss open problems.

**2. The regular lattice.** Throughout, for positive sequences $(a_m)$, $(b_m)$, we write $a_m \asymp b_m$ if the ratio $a_m/b_m$ is bounded away from zero and infinity. Also, we occasionally drop subscripts to lighten the notation, wherever there is no ambiguity.

2.1. *Bayesian detection.* We assume the uniform distribution over all paths, denoted by $\pi$. Equivalently, the distribution of the unknown path is that of an oriented symmetric random walk. We write $\mathbf{P}_\pi(\cdot) = \mathbf{E}_\pi \mathbf{P}_{1,p}(\cdot)$. As is well known, the test minimizing the risk (1.3) is the Neyman–Pearson test which rejects the null if and only if the likelihood ratio $L_m(X) = d\mathbf{P}_\pi(X)/d\mathbf{P}_0(X)$ exceeds 1 (the subscript $m$ refers here to the size of the problem). Here, the likelihood ratio is given by

$$(2.1) \qquad L_m(X) = 2^{-(m-1)} \sum_{p \in \mathcal{P}_m} e^{\mu X_p - m\mu^2/2},$$

where $X_p$ is defined in (1.4). Although $L_m(X)$ is an average over an exponentially large number of paths so that, at first sight, calculating this quantity may seem practically impossible, there is a recurrence relation which actually gives an algorithm for computing the likelihood ratio in a number of operations which is proportional to the number of nodes; see Section 6 for details. Note that the likelihood ratio $L_m(X)$ is closely related to the partition function of models of random polymers; see [11].

2.1.1. *Proof of Theorem 1.2: upper bound.* Assume $\mu_m m^{1/4} \to \infty$. This implies the existence of a sequence of real numbers $(h_m)$ tending to infinity and such that $\mu_m m^{1/4} h_m^{-1/2} \to \infty$. Define $\mathcal{S}(h_m)$ as the set of nodes obeying

$$\mathcal{S}(h_m) = \{(i,j) \in V_m : |j| \leq h_m \sqrt{m}\}.$$

In other words, $\mathcal{S}(h_m)$ is the intersection of $V_m$ with a strip of width $h_m\sqrt{m}$. Define $T_m$, the sum of the variables in the strip,

$$(2.2) \qquad T_m = \sum_{(i,j)\in\mathcal{S}(h_m)} X_{i,j},$$

and consider the test rejecting for appropriate large values of $T_m$ (determined below). With the assumption that $h_m$ is going to infinity, the oriented symmetric random walk $(i, S_i)_{0\leq i\leq m-1}$ is contained in $\mathcal{S}(h_m)$ with probability approaching 1. That is, if we define the event

$$A_m = \left\{\max_{0\leq i\leq m-1}|S_i| < h_m\sqrt{m}\right\},$$

then

$$(2.3) \qquad \lim_{m\to\infty} \mathbf{P}(A_m^c) = 0.$$

To see this, use Doob inequality for martingales to get

$$\mathbf{P}\left(\max_{0\leq i\leq m-1}|S_i| \geq h_m\sqrt{m}\right) \leq 2\frac{\mathbf{E}|S_{m-1}|}{h_m\sqrt{m}} \leq 2\frac{\sqrt{\mathbf{E}|S_{m-1}|^2}}{h_m\sqrt{m}} \leq \frac{2}{h_m}.$$

Let $n_m$ be the number of nodes in $\mathcal{S}(h_m)$ and note that $n_m = h_m m^{3/2}(1+o(1))$. Under $H_0$,

$$T_m \sim N(0, n_m),$$

while under $H_1$, conditionally on the event $A_m$, we have

$$T_m \sim N(m\mu_m, n_m).$$

It then follows from (2.3) and the fact that

$$\mu_m m n_m^{-1/2} \asymp \mu_m m^{1/4} h_m^{-1/2} \to \infty \qquad \text{as } m\to\infty$$

that the test with rejection region $|T_m| > m\mu_m/2$ obeys

$$\lim_{m\to\infty} \mathbf{P}_0(\text{Type I}) = 0 \quad \text{and} \quad \lim_{m\to\infty} \mathbf{P}_\pi(\text{Type II}) = 0.$$

That is, the test based on $T_m$ is asymptotically powerful. This completes the proof of the first part of Theorem 1.2.

2.1.2. *Proof of Theorem 1.2: lower bound.* It suffices to bound the Bayes risk from below. Note that

$$(2.4) \qquad B_m(\pi) := \inf_{\text{all tests}} \gamma_\pi(T_m) = \mathbf{P}_0(L_m \geq 1) + \mathbf{P}_\pi(L_m < 1),$$

where $L_m$ is the Bayes test (or likelihood ratio) $L_m(X) = (d\mathbf{P}_\pi/d\mathbf{P}_0)(X)$ and $E_0 L_m = 1$. A standard calculation shows that

$$(2.5) \qquad B_m(\pi) = 1 - \frac{\mathbf{E}_0|L_m - 1|}{2} \geq 1 - \frac{\sqrt{\mathbf{E}_0(L_m - 1)^2}}{2}.$$

Therefore, to show that the two hypotheses are asymptotically indistinguishable, it is sufficient to establish that, under the null, the variance of the likelihood ratio tends to zero.

Another standard calculation shows that the variance of $L_m$ is given by

$$(2.6) \qquad \mathbf{E}_0(L_m - 1)^2 = \mathbf{E}_0 L_m^2 - 1 = \mathbf{E} e^{\mu_m^2 N_m} - 1,$$

where $N_m$ is the number of crossings of two independent paths of length $m$ drawn from the prior. Hence, to derive a lower bound with this strategy, one needs to understand for which sequences $(t_m)$

$$(2.7) \qquad \lim_{m \to \infty} M_m(t_m) = 1,$$

where $M_m(t) := \mathbf{E} e^{t N_m}$, $t \in \mathbb{R}$, is the moment generating function of $N_m$.

Here, the prior is the distribution of a symmetric random walk and the reader may know that $\mathbf{E} N_m \asymp m^{1/2}$. Since

$$\mathbf{E} e^{t_m N_m} \geq 1 + t_m \mathbf{E} N_m,$$

this shows that it is necessary for the bound to be effective, to have $t_m m^{1/2} \to 0$ or, equivalently, $\mu_m m^{1/4} \to 0$. This is the correct asymptotic behavior, as we shall see next.

Let $(S_i)_{1 \leq i \leq m}$ and $(S'_i)_{1 \leq i \leq m}$ be two independent symmetric random walks (note the slight change of the range of indices which is of no consequence whatsoever). Observe that $\{S_i = S'_i\} = \{S_i - S'_i = 0\}$ so that we equivalently need to study the number $N_m$ of returns to zero of the difference process $(S_i - S'_i)_{1 \leq i \leq m}$, which is a Markov chain with the even integers as state space, and with jump probabilities to each neighbor equal to $1/4$ and probability to stay put equal to $1/2$. Therefore, the joint law of the difference process is that of $(S_{2i})_{1 \leq i \leq m}$, where, again, $S$ is a symmetric random walk (note the doubling of the interval together with the sampling at even times only). An immediate consequence is that

$$\mathbf{P}(N_m = k) = \mathbf{P}(|\{1 \leq i \leq m : S_{2i} = 0\}| = k).$$

The number of returns of a random walk to the origin has been well studied and we have from [14] and [15], Page 96 that

$$(2.8) \qquad \mathbf{P}(N_m = k) = \frac{1}{2^{2m-k}} \binom{2m-k}{m}.$$

The idea is now to develop a useful upper bound on the right side of (2.8) in order to estimate the moment generating function of $N_m$.

First, recall the classical refinement of the Stirling approximation to $n!$ (see [15], pages 50–53), which states that

$$\sqrt{2\pi} n^{n+1/2} e^{-n+1/(12n+1)} < n! < \sqrt{2\pi} n^{n+1/2} e^{-n+1/(12n)}.$$

Substituting this approximation into (2.8) (when expanding the binomial coefficients) yields

$$\mathbf{P}(N_m = k) \leq \frac{1}{\sqrt{\pi m}} \frac{(1 - k/2m)^{2m-k+1/2}}{(1 - k/m)^{m-k+1/2}}$$

(2.9)

$$= \frac{1}{\sqrt{\pi m}} \sqrt{\frac{1 - k/2m}{1 - k/m}} e^{-mg(k/m)},$$

where

$$g(t) = (1 - t)\log(1 - t) - 2(1 - t/2)\log(1 - t/2), \qquad 0 \leq t \leq 1.$$

For $t \in (0, 1)$, it holds that $d^2/dt^2(g(t) - t^2/4) > 0$ and, by convexity, the function $g(t) - t^2/4$ is above its tangent at the origin. This tangent is the line $y = 0$ since $g(0) = g'(0) = 0$, whence

$$g(t) \geq t^2/4 \qquad \forall t \in [0, 1).$$

Also, observe that $(1 - t/2)/(1 - t) \leq 1 + t$ for each $t \in [0, 1/2]$. Now, fix $0 < \epsilon < 1/2$. For $k \leq \epsilon m$, we have $(1 - k/2m)/(1 - k/m) \leq 1 + \epsilon$, while for $k < m$, one always has $\sqrt{\frac{1-k/2m}{m(1-k/m)}} \leq 1$. We then conclude that

(2.10) $$\mathbf{P}(N_m = k) \leq \begin{cases} \sqrt{\frac{1+\epsilon}{\pi m}} e^{-k^2/4m}, & k \leq \epsilon m, \\ \frac{1}{\sqrt{\pi}} e^{-k^2/4m}, & \epsilon m < k \leq m. \end{cases}$$

[The case $k = m$ in the above estimate is checked directly rather than from (2.9).]

The estimate (2.10) gives an upper bound on the moment generating function at $t_m$ since

$$M_m(t_m) \leq \sum_{k=0}^{\lfloor \epsilon m \rfloor} e^{t_m k} \frac{\sqrt{1+\epsilon}}{\sqrt{\pi m}} e^{-k^2/4m} + \sum_{k=\lfloor \epsilon m \rfloor + 1}^{m} e^{t_m k} \frac{1}{\sqrt{\pi}} e^{-k^2/4m}.$$

It is clear that if $t_m \to 0$ as $m \to \infty$, then the second term of the right side goes to zero as $m \to \infty$ so we focus on the first term. Using the monotonicity in $k$ of both $e^{t_m k}$ and $e^{-k^2 m/4}$, we have

$$\sum_{k=0}^{\lfloor \epsilon m \rfloor} e^{t_m k} \frac{1}{\sqrt{\pi m}} e^{-k^2/4m} \leq \frac{1}{\sqrt{\pi m}} + \sqrt{\frac{m}{\pi}} e^{t_m} \int_0^{\epsilon} e^{m t_m u} e^{-m u^2/4} \, du$$

$$= \frac{1}{\sqrt{\pi m}} + 2 e^{t_m} \int_0^{\epsilon \sqrt{m/2}} e^{\sqrt{2m} t_m u} \frac{e^{-u^2/2}}{\sqrt{2\pi}} \, du$$

$$\leq \frac{1}{\sqrt{\pi m}} + 2 e^{m t_m^2 + t_m} \mathbf{P}(Z > -\sqrt{2m} t_m),$$

where $Z$ is a standard normal random variable. It follows that if $t_m$ is chosen such that $\sqrt{m}t_m \to 0$ as $m \to \infty$, then

$$\lim_{m\to\infty} 2e^{mt_m^2 + t_m} \mathbf{P}(Z > -\sqrt{2m}t_m) = 1$$

and thus $\lim_{m\to\infty} M_m(t_m) = 1$. In conclusion, we have proven that

$$(2.11) \qquad \lim_{m\to\infty} \mu_m m^{1/4} = 0 \implies \liminf_{m\to\infty} B_m(\pi) \geq 1.$$

This proves the second part of Theorem 1.2.

2.2. *Minimax detection.* Just as in the Bayesian case, we first prove the upper bound by constructing a test which allows us to detect reliably when $\mu_m$ decays slower than $(\log m)^{-1/2}$, and then study the lower minimax bound.

2.2.1. *Proof of Theorem 1.1: upper bound.* Consider a simple test statistic of the form

$$(2.12) \qquad T_m = \sum_{(i,j)\in V_m} w_{i,j} X_{i,j}, \qquad w_{i,j} := w_i = \frac{\lambda_m}{i+1}.$$

Hence, $T_m$ is a weighted sum of the values at the vertices of the graph. For convenience, we fix $\lambda_m$ so that $\sum_{0 \leq i \leq m-1} w_i = 1$. Note that $\lambda_m = (\log m)^{-1}(1 + o(1))$. Under $H_1$, the mean of $T_m$ is given by $\mu_m \sum_{0 \leq i \leq m-1} w_i = \mu_m$ and since the $X_{i,j}$'s have identical variance under both $H_0$ and $H_1$, we have

$$\mathrm{Var}_0(T_m) = \mathrm{Var}_{1,p}(T_m) = \sum_{(i,j)\in V_m} w_{i,j}^2$$

$$= \sum_{0\leq i\leq m-1} (i+1)w_i^2 = \sum_{0\leq i\leq m-1} \frac{\lambda_m^2}{i+1} = \lambda_m.$$

Hence,

$$T_m \sim_{H_0} N(0, \lambda_m) \quad \text{and} \quad T_m \sim_{H_{1,p}} N(\mu_m, \lambda_m),$$

under any alternative. Consider the test which rejects the null whenever $T_m > \mu_m/2$. The risk of this test is then equal to

$$\gamma(T_m) = 2P(N(0,1) > \tfrac{1}{2}\mu_m \lambda_m^{-1/2}) \implies \lim_{m\to\infty} \gamma(T_m) = 0$$

when $\mu_m \lambda_m^{-1/2} \to \infty$ or, equivalently, when $\mu_m \sqrt{\log m} \to \infty$. This proves the first part of Theorem 1.1.

2.2.2. *Proof of Theorem* 1.1: *lower bound.* The idea for obtaining a lower bound is to exhibit a prior on $H_1$ which makes the Bayesian detection problem as hard as possible. Consider a prior $\pi$ on $H_1$ (here a distribution on the set of paths). Then, for all tests $T_m$,

$$\gamma(T_m) \geq B_m(\pi),$$

where $B_m(\pi)$ is the risk of the Bayes test,

$$B_m(\pi) = \mathbf{P}_0(L_m \geq 1) + \mathbf{P}_\pi(L_m < 1).$$

Our strategy is to construct a prior on the family of paths with a low predictability profile, that is, a process whose location in the future is hard to predict from its current state and history.

*The predictability profile of a stochastic process.* The concept of the predictability profile was first introduced in [7].

DEFINITION 2.1. The predictability profile of a stochastic process $(S_n)_{n \geq 1}$ is defined by

$$(2.13) \qquad \mathrm{PRE}_S(k) = \sup \mathbf{P}(S_{n+k} = x | S_0, \ldots, S_n),$$

where the supremum is taken over all positions and histories.

We will consider nearest-neighbor walks which are defined as processes with increments equal to $\pm 1$. Improving upon earlier results of Benjamini, Pemantle and Peres [7], Häggström and Mossel [17], Theorem 1.4, proved the following.

THEOREM 2.2. *Suppose $(f_k)_{k \geq 1}$ is a decreasing positive sequence such that $\sum_{k \geq 1} f_k/k < \infty$. There then exists a nearest-neighbor process starting at $S_0 = 0$ and obeying*

$$(2.14) \qquad \mathrm{PRE}_S(k) \leq \frac{C}{k f_k}$$

*for all $k \geq 1$ and some positive constant $C$.*

C. Hoffman proved in [18] that this is sharp in the sense that if $(f_k)$ is a decreasing positive sequence with $\sum_{k \geq 1} f_k/k = \infty$, then the predictability profile (2.14) is impossible to achieve.

In what follows, we will need a quantitative, finite version of Theorem 2.2. This is achieved by using a concrete prior, introduced in [17], which gives the predictability profile below.

LEMMA 2.3 ([17], Proposition 3.1). *Fix a sequence $(a_j)_{j \geq 0}$ obeying $\sum_{j \geq 0} a_j < 1$. There then exists a nearest-neighbor process $(S_n)_{n \geq 0}$ obeying*

$$(2.15) \qquad \mathrm{PRE}_S(k) \leq \frac{20}{k a_{\lfloor \log_2(k/2) \rfloor}} \qquad \text{for all } k = 1, 2, \ldots.$$

The construction of the process and the proof of (2.15) may be found in the Appendix. Later, we will consider a prior on paths obeying (2.15) for suitable values of the sequence $(a_j)_{j \geq 0}$.

*Predictability profiles and numbers of intersections.* Hereafter, we consider stochastic processes with a finite horizon, that is, $(S_i)_{0 \leq i \leq m-1}$. In the sequel, we will need to estimate the number of times two independent processes drawn from a prior with prescribed predictability profile cross each other. From the proof of [7], Lemma 3.1, we state the following

LEMMA 2.4. *Let $B$ be such that*

$$(2.16) \qquad \sum_{1 \leq k \leq \lfloor m/B \rfloor} \mathrm{PRE}_S(kB) \leq \theta < 1.$$

*Then, for any sequence $(v_n)_{0 \leq n \leq m-1}$ and all $k \geq 1$, the distribution of the total number of intersections between $(S_n)$ and $(v_n)$ obeys*

$$(2.17) \qquad \mathbf{P}(|S \cap v| \geq k) \leq B \cdot \theta^{k/B}, \qquad |S \cap v| := |\{n : S_n = v_n\}|.$$

We emphasize that the lemma is valid even if the sequence $(v_n)_{n \geq 0}$ does not determine a nearest-neighbor path.

We now prove the lower bound in Theorem 1.1 by providing a lower bound for the Bayes risk $B_m(\pi)$ for the prior $\pi$ given by Lemma 2.3, and with the sequence

$$(2.18) \qquad a_j = a_j(m) := \begin{cases} 1/(3 \log_2 m), & j \leq \log_2 m, \\ 0, & j > \log_2 m. \end{cases}$$

With the above choice, $\sum_{j \geq 0} a_j \leq \frac{\log_2 m + 1}{3 \log_2 m} < 1/2$ for $m > 4$.

As in the analysis of the Bayes risk [see (2.5),(2.6)], we employ the simple bound

$$(2.19) \quad B_m(\pi) \geq 1 - \frac{\sqrt{\mathbf{E}_0(L_m-1)^2}}{2}, \qquad \mathbf{E}_0(L_m-1)^2 = \mathbf{E} e^{\mu_m^2 N_m} - 1,$$

where $L_m$ is the likelihood ratio and $N_m$ is the number of crossings of two independent paths drawn from the prior $\pi$. We compute

$$\sum_{k \geq 1} e^{\mu_m^2 k} \mathbf{P}(N_m = k) = \sum_{1 \leq k \leq K-1} e^{\mu_m^2 k} \mathbf{P}(N_m = k)$$
$$+ \sum_{k \geq K} e^{\mu_m^2 k} [\mathbf{P}(N_m \geq k) - \mathbf{P}(N_m \geq k+1)]$$

and, summing by parts, deduce that

$$\mathbf{E}e^{\mu_m^2 N_m} \leq e^{\mu_m^2(K-1)} + [1 - e^{-\mu_m^2}] \sum_{k \geq K} \mathbf{P}(N_m \geq k)e^{\mu_m^2 k}.$$

With the choice (2.18), Lemma 2.3 gives

$$\mathrm{PRE}_S(k) \leq (60 \log_2 m)/k, \qquad k = 1, 2, \ldots.$$

In particular, with $B = B_m = 120(\log m)^2/\log 2$, we have

$$\sum_{k=1}^{\lfloor m/B_m \rfloor} \mathrm{PRE}_S(kB_m) \leq \tfrac{1}{2}.$$

Applying Lemma 2.4 yields

$$\mathbf{E}_0 L_m^2 \leq e^{\mu_m^2(K-1)} + [1 - e^{-\mu_m^2}] B_m \sum_{k \geq K} e^{\mu_m^2 k} 2^{-k/B_m}$$

$$\leq e^{\mu_m^2(K-1)} + [1 - e^{-\mu_m^2}] B_m \frac{a_m^K}{1 - a_m}, \qquad a_m = e^{\mu_m^2} 2^{-1/B_m} < 1,$$

where the last inequality is due to the fact that $\lim_{m \to \infty} \mu_m^2 B_m = 0$ [since $\mu_m (\log m)(\log \log m)^{1/2} = o(1)$]. Further,

$$\liminf_{m \to \infty}(-B_m \log a_m) = \log 2 \implies \frac{1}{1 - a_m} \leq \frac{1}{1 - e^{-\log 2/(2B_m)}} \leq c_1 B_m$$

for some constant $c_1$ and all $m$ large. It follows that for some constant $c_2$ and all $m$ large,

$$\mathbf{E}_0 L_m^2 \leq e^{\mu_m^2 K} + c_2 \mu_m^2 B_m^2 e^{-K(\log 2)/2B_m}.$$

Taking $K = K_m = 2(B_m \log B_m)/\log 2$ yields, for some constant $c_3$,

$$\mathbf{E}_0 L_m^2 \leq e^{c_3 \mu_m^2 B_m \log B_m} + O(\mu_m^2 B_m) \to 1 \qquad \text{as } m \to \infty.$$

Together with (2.19), this concludes the proof of Theorem 1.2.

**3. The complete binary tree.** In this section, we prove Theorem 1.3. For the upper bound, we show that the GLRT is asymptotically powerful if $\mu_m = \mu > \sqrt{2 \log 2}$ and that a closely related test is asymptotically powerful if $\mu_m = \mu = \sqrt{2 \log 2}$. For the lower bound, we study the likelihood ratio under the uniform prior on paths using a martingale approach.

We start by considering the GLRT, which is based on $M_m = \max\{X_p : p \in \mathcal{P}_m\}$, $X_p$ being defined in (1.4). We first show that under the null hypothesis, the GLRT obeys

$$\mathbf{P}_0(M_m \geq m\sqrt{2 \log 2}) \to 0, \qquad m \to \infty.$$

This is, in fact, a simple application of Boole's inequality and a standard bound on the tail of the normal distribution

$$\mathbf{P}(N(0,1) > t) \leq \frac{1}{\sqrt{2\pi}} \frac{e^{-t^2/2}}{t}.$$

Indeed,

$$\mathbf{P}_0(M_m \geq m\sqrt{2\log 2}) \leq 2^{m-1} \mathbf{P}_0(X_p \geq m\sqrt{2\log 2})$$

(3.1)

$$= 2^{m-1} \mathbf{P}(N(0,1) \geq \sqrt{2m\log 2}) \leq \frac{1}{4\sqrt{\pi m \log 2}}.$$

In fact, $M_m/m \to \sqrt{2\log 2}$ a.s.; see [23], Section 3. Under any alternative $\mathbf{P}_{1,p}$ with $\mu > \sqrt{2\log 2}$, however, the GLRT obeys

(3.2) $$\mathbf{P}_{1,p}(M_m > m\sqrt{2\log 2}) \to 1, \qquad m \to \infty.$$

Indeed, if $p$ is the path along which the mean is elevated, $M_m \geq X_p$ and $X_p/m$ is normally distributed with mean $\mu$ and variance $1/m$.

If $\mu = \sqrt{2\log 2}$, the same argument gives

$$\liminf_{m\to\infty} \mathbf{P}_{1,p}(M_m > m\sqrt{2\log 2}) \geq \tfrac{1}{2}$$

for each path $p$ instead of (3.2). This is not quite enough to conclude that $H_0$ and $H_1$ can be separated with probability approaching 1. However, taking $m_k = 2^k$, from (3.1) and Borel–Cantelli, we have that

$$\mathbf{P}_0(M_{m_k} \geq m_k\sqrt{2\log 2} \text{ infinitely often}) = 0,$$

while standard estimates for random walks imply that

$$\mathbf{P}_{1,p}(M_{m_k} \geq m_k\sqrt{2\log 2} \text{ infinitely often}) = 1$$

for each $p$ and [because the increments $X_p(m_k) - X_p(m_{k-1})$ are exponentially mixing] even

$$\liminf_{k\to\infty} \frac{1}{k} \sum_{i=1}^{k} \mathbf{1}_{\{M_{m_i} \geq m_i\sqrt{2\log 2}\}} \geq \frac{1}{2}, \qquad \mathbf{P}_{1,p}\text{-a.s.}$$

Therefore, the test which computes, along the sequence $m_k$, the number of times $M_{m_k} \geq m_k\sqrt{2\log 2}$, declaring $H_0$ if this number is less than $k/4$ and $H_1$ otherwise, has asymptotic full power.

In conclusion, the GLRT (or its variant) has asymptotic full power if $\mu \geq \sqrt{2\log 2}$.

We now turn to studying the likelihood ratio under the uniform prior $\pi$ on paths

$$L_m = 2^{-(m-1)} \sum_{\text{all paths } p} e^{\mu X_p - m\mu^2/2}$$

and show that for $\mu < \sqrt{2\log 2}$, its risk

$$B_m(\pi) = \mathbf{P}_0(L_m \geq 1) + \mathbf{P}_\pi(L_m < 1)$$

is bounded away from 0. A lower bound such as (2.5) would not suffice here since we want to recover the same threshold $\sqrt{2\log 2}$. Instead, we turn to martingale methods. Such methods have been used for years (see, e.g., [8, 12]). Here, we follow the presentation found in [9].

A simple calculation shows that

$$B_m(\pi) = 1 - \mathbf{E}_0(1 - L_m)_+.$$

Let $|v|$ denote the distance of a vertex $v$ from the root. By Proposition 1 in [9], we know that under $H_0$, $L_m$ is a nonnegative martingale with respect to the filtration $\mathcal{F}(X_v : |v| \leq m-1)$, which converges pointwise to a finite, nonnegative random variable $L_\infty$. Hence, by dominated convergence,

$$\lim_{m\to\infty} B_m(\pi) = 1 - \mathbf{E}_0(1 - L_\infty)_+.$$

Applying Proposition 2 in [9], we have that for $\mu < \sqrt{2\log 2}$, $L_m$ is uniformly integrable and, therefore, $\mathbf{E}_0 L_\infty = 1$. Hence, $\mathbf{P}_0(L_\infty = 0) < 1$ and, consequently,

$$\lim_{m\to\infty} B_m(\pi) > 0.$$

Finally, we briefly argue that if $\mu_m \to 0$, then every sequence of tests is asymptotically powerless. Here, it is enough to use the bound (2.5). It therefore suffices to prove that $\mathrm{Var}_0(L_m) \to 0$ as $m \to \infty$. Just as in (2.19), $\mathrm{Var}_0(L_m) = \mathbf{E}_\pi e^{\mu_m^2 N_m} - 1$, where $N_m$ is the number of crossings between two random paths drawn from the prior $\pi$. Here, $\mathbf{P}(N_m = k) = 2^{-k}$, $1 \leq k \leq m-1$, and $\mathbf{P}(N_m = m) = 2^{-m+1}$. In short, the distribution of $N_m$ is that of a truncated geometric random variable with probability of success equal to $1/2$. Set $\tau_m = e^{\mu_m^2}/2$, which is less than 1 for $m$ large. We compute

$$\mathrm{Var}_0(L_m) = \frac{(2\tau_m - 1)(1 - \tau_m^m)}{1 - \tau_m} \leq \frac{2\tau_m - 1}{1 - \tau_m}.$$

It is now clear that $\mathrm{Var}_0(L_m) \to 0$ when $\tau_m \to 1/2$ or, equivalently, when $\mu_m \to 0$. This completes the proof of the theorem.

**4. Extension to exponential families.** While the previous sections studied the detection problem assuming a Gaussian distribution at the nodes of the graph, it is now time to emphasize that our results hold more generally. In fact, one can obtain similar conclusions for exponential models as well.

Letting $F_0$ be a distribution on the real line, we define $F_\theta$ as the exponential family with associated density $\exp(\theta x - \log\varphi(\theta))$ with respect to $F_0$.

Note that by definition, $\varphi(\theta) = \mathbf{E}_{F_0}[\exp(\theta X)]$, where $\mathbf{E}_{F_0}$ is the expectation under the distribution $F_0$. We always assume that $\varphi(\theta) < \infty$ for $\theta$ in a neighborhood of 0; further restrictions are mentioned when needed.

Under the null hypothesis, we assume that all the nodes are i.i.d. with distribution $F_0$, while under $H_{1,m}$, there is a path along which the nodes are i.i.d. with distribution $F_{\theta_m}$, $\theta_m > 0$, and distribution $F_0$ away from the path. The question is, of course, for what values of $\theta_m$ one can reliably detect this path. To connect this general set-up with the previously studied special case, set $\psi(\theta) = \log \varphi(\theta)$ and recall that

$$\mu(\theta) := \mathbf{E}_{F_\theta} X = \psi'(\theta) \quad \text{and} \quad \sigma^2(\theta) := \operatorname{Var}_{F_\theta} X = \psi''(\theta).$$

With this notation, the mean shift is equal to

$$\mu(\theta) - \mu(0) = \psi'(\theta) - \psi'(0) = \psi''(0)(\theta + o(\theta)).$$

In other words, the value of a small mean shift is just about proportional to $\theta$. [In the Gaussian case, $\mu(\theta) = \theta$ and $\log \varphi(\theta) = \theta^2/2$.]

4.1. *The regular lattice with an exponential family at the nodes.* We first consider the minimax detection problem, and extend Theorem 1.1.

THEOREM 4.1. *Suppose that $\theta_m \sqrt{\log m} \to \infty$ as $m \to \infty$. There is then a sequence of tests which is asymptotically powerful. Conversely, suppose that $\theta_m \log m \sqrt{\log \log m} \to 0$ as $m \to 0$. Then, every sequence of tests $(T_m)$ is asymptotically powerless.*

In summary, one can reliably detect a path as long as the mean shift $\mu(\theta_m) - \mu(0) \gg (\log m)^{-1/2}$, while this is impossible if—ignoring the $\sqrt{\log \log m}$ factor—$\mu(\theta_m) - \mu(0) \ll (\log m)^{-1}$.

As an example, consider the case where we have exponentially distributed random variables; under the null, the node variables are exponentially distributed with mean 1, while under the alternative hypothesis, there is a path along which the node variables are exponentially distributed with mean $1 + \mu_m$. Let $F_0$ be the density of the exponential with mean 1. The density of an exponential random variable with mean $1 + \mu$ with respect to $F_0$ is given by

$$(1+\mu)^{-1} \exp(\mu x/(1+\mu)) := \exp(\theta x - \log \varphi(\theta)),$$

with

$$\theta = \frac{\mu}{1+\mu}, \qquad \varphi(\theta) = \frac{1}{1-\theta}.$$

For this exponential model, one can reliably detect a mean shift $\mu_m$ if it is significantly larger than $(\log m)^{-1/2}$, while this is impossible if it is much smaller than $(\log m)^{-1}$.

PROOF OF THEOREM 4.1. The proof is similar to that of Theorem 1.1. For the upper bound, we consider the same statistic (2.12) as before, $T_m := \sum_{(i,j) \in V_m} w_{i,j} X_{i,j}$, with exactly the same choice of weights. First, observe that for any path $p$ from $H_1$, the mean difference obeys

$$\mathbf{E}_{1,p}(T_m) - \mathbf{E}_0(T_m) = \mu(\theta_m) - \mu(0).$$

As for the variances, we have

$$\mathrm{Var}_0(T_m) = \sigma^2(0) \sum_{(i,j) \in V_m} w_{i,j}^2 = \lambda_m \sigma^2(0)$$

and for any alternative in $H_1$,

$$\mathrm{Var}_{1,p}(T_m) = \sigma^2(0) \sum_{(i,j) \in V_m} w_{i,j}^2 + [\sigma^2(\theta_m) - \sigma^2(0)] \sum_{0 \le i \le m-1} w_i^2$$

$$= \lambda_m \sigma^2(0) + [\sigma^2(\theta_m) - \sigma^2(0)] O(\lambda_m^2).$$

Recall that $\lambda_m = (\log m)^{-1}(1 + o(1))$. Using Chebychev's inequality, we see that the probabilities of Type I and Type II errors go to zero as soon as $[\mu(\theta_m) - \mu(0)]\lambda_m^{1/2} \to \infty$ as $m \to \infty$. The first part of the theorem follows from $\mu(\theta_m) - \mu(0) = \theta_m \mathrm{Var}_{F_0}(X)(1 + o(1))$. That is, if the mean shift time $\sqrt{\log m}$ increases to infinity, then the probability of each type of error goes to zero.

For the lower bound, we consider the same prior distribution on the family of paths. For exponential models, the variance of the likelihood ratio $L_m$ is given by

(4.1) $$\mathrm{Var}_0(L_m) = \mathbf{E}[\lambda(\theta_m)^{N_m}] - 1, \qquad \lambda(\theta) = \frac{\varphi(2\theta)}{\varphi(\theta)^2} > 1,$$

where, again, $N_m$ is the number of crossings of two independent paths drawn from the prior, or

$$\mathrm{Var}_0(L_m) = \mathbf{E} e^{\alpha^2(\theta_m) N_m} - 1 \qquad \alpha(\theta) = \sqrt{\log \lambda(\theta)}.$$

This is the same expression as before [cf. (2.19)] and our previous analysis shows the existence of a prior with the property

$$\lim_{m \to \infty} \alpha(\theta_m) \log m \sqrt{\log \log m} = 0 \implies \lim_{m \to \infty} \mathrm{Var}_0(L_m) = 0,$$

which implies that the Bayes test is asymptotically powerless. It is now not difficult to see that for exponential models, $\lambda(\theta) = 1 + O(|\theta|^2)$ so that $\alpha(\theta) = O(\theta)$ for $\theta$ close to zero. As a consequence,

$$\lim_{m \to \infty} \theta_m \log m \sqrt{\log \log m} = 0 \implies \lim_{m \to \infty} \alpha(\theta_m) \log m \sqrt{\log \log m} = 0,$$

which establishes the second part of the theorem. □

Not surprisingly, the same extension also holds in the Bayesian set-up.

THEOREM 4.2. *Consider the uniform prior on paths. Suppose that $\theta_m m^{1/4} \to \infty$ as $m \to \infty$. The Bayes test is then asymptotically powerful. Conversely, if $\theta_m m^{1/4} \to 0$ as $m \to 0$, then the Bayes risk tends to 1 and every sequence of tests $(T_m)$ is asymptotically powerless.*

The proof follows that of Theorems 1.2 and 4.1. We omit the details.

4.2. *The tree with an exponential family at the nodes.* Following [9], define the function $f$ as

$$f(\theta) = \frac{1}{\theta} \log(2\varphi(\theta)). \tag{4.2}$$

By Lemma 4 in [9], $f$ either attains its unique minimum or $f$ is strictly decreasing on $(0, \infty)$. In any case, we denote by $\theta^\star \in (0, \infty]$ the value where $f$ is minimum.

THEOREM 4.3. *Assume that $\varphi(\theta) < \infty$ in a neighborhood of $\theta^\star$. If $\theta_m = \theta > \theta^\star$, then the GLRT is asymptotically powerful. If $\theta_m = \theta < \theta^\star$, then there does not exist any asymptotically powerful sequence of tests. If $\theta_m \to 0$, then all sequences of tests are powerless. Finally, if $\theta_m = \theta^\star$, then there exists a sequence of asymptotically powerful tests.*

For exponential random variables, $\varphi(\theta) = 1/(1-\theta)$ and we numerically compute $\theta^\star \approx .63$. In terms of mean shift (see above), we have $\mu(\theta^\star) - \mu(0) = 1/(1-\theta^\star) - 1 \approx 1.70$. The mean difference along the unknown path must exceed approximately 1.70 to be reliably detectable.

For Bernoulli random variables, $F_\theta = \text{Bernoulli}(e^\theta/(1+e^\theta))$, the function $f$ is decreasing on $(0, \infty)$ and, therefore, $\theta^\star = \infty$. Theorem 4.3 then implies that no asymptotically powerful sequences of tests exist for testing fair coin tossing at the nodes versus biased coin tossing with parameter $q \in (1/2, 1)$ along a path. Note that the situation drastically changes when $q = 1$: in this case, the nodes with value 1 that are connected to the root node through a path of nodes of value 1 form a critical branching process (with an expected number of descendants at each node equal to 1) which, therefore, eventually dies out. Under $H_1$, however, there is always a path of length $m$ starting from the origin and with all 1's. Hence, the test that declares $H_1$ if one finds such a path and $H_0$ otherwise is asymptotically powerful.

PROOF OF THEOREM 4.3. The proof is very similar to that of Theorem 1.3. We start with the upper bound, assuming $\theta^\star < \infty$. Define $\xi(t) = \inf_{\theta > 0} \varphi(\theta) e^{-t\theta}$. Note that

$$\begin{aligned} \xi(t) = 1/2 &\iff \inf_{\theta > 0}(\log(2\varphi(\theta)) - \theta t) = 0 \\ &\iff t = \inf_{\theta > 0} f(\theta) = f(\theta^\star). \end{aligned} \tag{4.3}$$

Because $\varphi(\theta) < \infty$ in a neighborhood of $\theta^\star$, we can replace the estimate (3.1) by the Bahadur–Rao bound [4], which yields

$$\mathbf{P}_0(M_m \geq m\xi^{-1}(1/2)) \leq 2^{m-1}\mathbf{P}_0(X_p \geq m\xi^{-1}(1/2)) \leq \frac{C}{\sqrt{m}}$$

for some constant $C$. (In fact, under our assumptions, $M_m/m \to \xi^{-1}(1/2)$ a.s., by the argument in [23], Section 3.) This estimate and (4.3) imply that

(4.4) $$\mathbf{P}_0(M_m \geq mf(\theta^\star)) \leq \frac{C}{\sqrt{m}}.$$

We now study the behavior of $M_m/m$ under $H_1$. Let $p$ be the path along which the nodes are sampled from the distribution $F_\theta$. The strong law of large numbers then shows that $\lim_{m\to\infty} X_p/m = E_{F_\theta}X$ a.s. and, therefore,

$$\liminf_{m\to\infty} \frac{M_m}{m} \geq \frac{d}{d\theta}(\log\varphi(\theta)) \qquad \text{a.s.}$$

The derivative obeys $d/d\theta(\log\varphi(\theta)) = \varphi'(\theta)/\varphi(\theta) > f(\theta^\star)$ if and only if $\theta > \theta^\star$. This equivalence follows from the identity

$$d/d\theta \log(\varphi(\theta)) - f(\theta^\star) = \theta f'(\theta) + f(\theta) - f(\theta^\star).$$

Since $f$ is decreasing on $(0,\theta^\star)$ and strictly increasing on $(\theta^\star,\infty)$, the right-hand side has the sign of $\theta - \theta^\star$. This analysis shows that the GLRT has asymptotic full power if $\theta > \theta^\star$, and the argument for handling $\theta = \theta^\star$ is the same as in the Gaussian case, using the full power of (4.4).

The study of the likelihood ratio under the uniform prior over paths is identical to that in the Gaussian case, with the exception that when proving the uniform integrability of the martingale $L_m$, we use Biggins's theorem (in the form given in [21]—noting the condition $\varphi(\theta) < \infty$ in a neighborhood of $\theta^\star$) instead of using Proposition 2 from [9]. [The latter proposition requires that $\varphi(\theta)$ be finite for all $\theta > 0$, or at least for $\theta = 2\theta^\star$.]  □

**5. Extension to other graphs.** This section emphasizes that results are available for other graphs and, in particular, for the analog of the regular lattice in higher dimensions.

- *Regular lattice in dimension $d' = d+1$.* This is the graph with vertex set

  $$V = \{(i, j_1, \ldots, j_d) : 0 \leq i, -i \leq j_k \leq i \text{ and } j_k \text{ has the parity of } i\}$$

  and oriented edges $(i, j_1, \ldots, j_d) \to (i+1, j_1 + s_1, \ldots, j_d + s_d)$, where $s_k = \pm 1$.

Consider a distribution from the exponential family at the nodes and the uniform prior on paths. In this case, the likelihood ratio has been studied in dimension $d+1$—under the name of the *partition function*—in the context of directed random polymers. Martingale methods work well in this context and the behavior of the likelihood ratio for $d \geq 3$ is similar to the behavior of the likelihood ratio for the tree that we studied in Section 3; see [11], Proposition 3.2.1. In particular, for $d \geq 3$, there are no asymptotically powerful sequences of tests if $\theta_m = \theta$ obeys $\lambda(\theta)\rho_d < 1$, where $\lambda(\theta)$ is defined as in (4.1) and $\rho_d$ is the return probability of a symmetric random walk in dimension $d$. (The results for $d = 2$ only imply that the Bayes risk tends to zero if $\theta_m = \theta > 0$.) In contrast, the minimax risk does not go to zero here and this follows from the construction of a prior with low predictability profile. We give a general statement in Theorem 5.3.

To establish a general result, we work with a connected graph (directed or undirected), with one vertex marked that we call the origin, and, as before, we let $\mathcal{P}$ be the set of self-avoiding paths starting at the origin and $\mathcal{P}_m \subset \mathcal{P}$ be the subset of paths of length $m$. Under the null hypothesis, all the nodes are i.i.d. $F_0$, while under the alternative, there is a path in $\mathcal{P}_m$ along which the nodes are i.i.d. $F_1$. We assume throughout that $F_1$ is absolutely continuous with respect to $F_0$; otherwise, the detection problem becomes trivial.

DEFINITION 5.1. A distribution $\pi$ on $\mathcal{P}$ is said to have an *exponential intersection tail with parameter* $\eta \in (0,1)$ if there exists $C > 0$ such that if $N$ is the number of crossings of two independent samples from $\pi$, then

$$\mathbf{P}(N \geq k) \leq C \cdot \eta^k \qquad \forall k \geq 1.$$

The regular lattice with $d \geq 2$ (i.e., $d' \geq 3$) admits a measure on paths with an exponential intersection tail [7], Theorem 1.3. Note that a summable predictability profile implies an exponential intersection tail.

DEFINITION 5.2. Let $L = dF_1/dF_0$ be the likelihood ratio for testing $F_1$ versus $F_0$ at a single node. The Pearson $\chi^2$-distance between $F_0$ and $F_1$ is defined as $\chi^2(F_0, F_1) = \text{Var}_0(L)$.

With these definitions, we have the following general statement.

THEOREM 5.3. *Suppose that there is a distribution $\pi$ on $\mathcal{P}$ having an exponential intersection tail with parameter $\eta$. Then, if $\chi^2(F_0, F_1) < \eta^{-1} - 1$, there are no asymptotically powerful sequences of tests.*

The proof does not require any argument that we have not already presented, and is omitted. For exponential variables, $\chi^2(F_0, F_\theta) = \lambda(\theta) - 1$, where $\lambda(\theta)$ is defined as in (4.1) and, therefore, no asymptotically powerful sequences of tests exist if $\lambda(\theta)\eta < 1$.

Theorem 5.3 provides a lower bound on the minimax threshold for reliable detection. For an upper bound, suppose, for example, that the variables are exponentially distributed and assume that $\#\mathcal{P}_m = O(\delta^m)$ for some positive constant $\delta$; for instance, $\delta = 2^d$ works for the regular lattice in dimension $d+1$. Application of Boole's inequality and the law of large numbers shows that under those assumptions, the GLRT is asymptotically powerful if $\xi(t)\delta > 1$, where, again, $\xi(t) = \inf_{\theta>0} \varphi(\theta) e^{-t\theta}$.

**6. Numerical experiments.** We now explore the empirical performance of some of the detection methods we proposed for the regular lattice. The variables at the nodes are independent Gaussians. To measure the performance, we fix the probability of Type I error at 5% and estimate the power or detection rate, that is, the probability of deciding in favor of the alternative $H_1$ when $H_1$ is true. This power function was estimated at values of the mean shift $\mu$ (the mean of the node variables along the path) at which this function is varying.

6.1. *Bayesian detection under the uniform prior.* We first consider detection under the uniform prior on paths. We compare the performance of the Bayes test, the GLRT and the test based on the strip statistic which was used in the proof of the upper bound in Theorem 1.2. The Bayes test is optimal in this setting and we recall that the strip statistic was shown to achieve the optimal detection rate. This paper did not theoretically analyze the performance of the GLRT in this situation, however, and we would like to do so empirically.

6.1.1. *Computing the Bayes statistic.* As emphasized earlier, there exists a rapid algorithm for calculating the Bayes statistic $L_m(X)$ [(2.1)]. Consider any node $v = (i, j)$ ($0 \le i \le m-1$ and $j$ has the parity of $i$) and let $\mathcal{P}^{\mathrm{End}}(v)$ be the set of paths starting at the root $(0,0)$ and ending at the node $v$. Set

$$Y(v) := 2^{-i} \sum_{p \in \mathcal{P}^{\mathrm{End}}(v)} e^{\mu X_p - (i+1)\mu^2/2}.$$

With this notation, $L_m(X)$ is the sum of $Y$ over all the terminal nodes $v$ for which $i = m-1$. Now, observe the recurrence

(6.1) $$Y(v) = e^{\mu X_v - \mu^2/2} \cdot \frac{Y(v^+) + Y(v^-)}{2},$$

where $(v^+, v^-)$ are the two predecessors of $v$ in the graph, that is, the two nodes from which one can reach $v$ in one step. [By convention, set $Y(v^\pm) = 0$ if $v^\pm$ is outside the grid.] This recurrence shows that one can compute the Bayes statistics in $O(m^2)$ flops.

For each value of $\mu$ and $m$, then, we simulated the Bayes statistic under $H_0$ and $H_1$ using 2,000 realizations for each. Here and below, each realization uses a new path realization drawn from the uniform distribution.

6.1.2. *Simulating the strip statistic.* For a positive integer $B$, the strip statistic $T_{m,B}$ is the sum of the random variables falling in the centered strip of length $m$ and width $2B+1$,

$$T_{m,B} = \sum_{0 \leq i \leq m-1} \sum_{j:|j| \leq \min(i,B)} X_{i,j}.$$

Under $H_0$, $T_{m,B} \sim N(0, n_{m,B})$ where $n_{m,B}$ is the number of vertices in the strip, while under $H_1$ (fixed path), $T_{m,B} \sim N(\mu \cdot R_{m,B}, n_{m,B})$, where $R_{m,B}$ is the number of vertices inside the strip that the random path visits. Therefore, one can simulate $T_{m,B}$ by taking one realization of $R_{m,B}$, multiplying it by $\mu$ and adding an independent mean-zero Gaussian variable.

It remains to choose the width of the strip. We ran simulations with $B = \nu\sqrt{m}$ for $\nu = 0.75, 1, 2, 3$. Among these values, $B = 2\sqrt{m}$ gave the best performance (at least for the graph sizes we considered). Finally, for a fixed $\mu$ and $m$, we used 5,000 realizations of the test statistic to estimate the detection rate.

6.1.3. *Simulating the GLRT.* The GLRT statistic rejects for large values of $M_m = \max\{X_p : p \in \mathcal{P}_m\}$. This statistic can be calculated rapidly using dynamic programming; for example, Dijkstra's algorithm [1] has here a computational complexity proportional to the number of nodes. For each graph size, the threshold corresponding to a Type I error probability approximately equal to .05 and the detection rate for a fixed $\mu$ were based on 10,000 and 1,000 realizations, respectively.

6.1.4. *Comparing the tests.* To compare the three tests, one can estimate the value of the mean shift which gives a detection rate of about 95% from graphs plotting the detection rates versus $\mu$ (see Figure 6). Call this quantity $\mu_{0.95}$. Table 1 shows $\mu_{0.95}$ for the Bayes test, the test based on the strip statistic test and the GLRT for different graph sizes. As expected, the Bayes test outperforms the other two, but one needs to recall that those tests do not require information about the parameter $\mu$, while the Bayes test does. Figure 3 shows a log-log plot of $\mu_{0.95}$ as a function of $m$, together with least-squares line fits. The slope of the line is $-0.255$ for the Bayes test and $-0.246$ for the strip test. Both of these values are quite close to the $-1/4$ exponent one finds in Theorem 1.2. For the GLRT, the slope is about $-0.16$. This suggests that the strip statistic test might eventually outperform the GLRT for sufficiently large graphs. The fitted lines meet at approximately $m = 2^{20} \approx 10^6$, but it would be computationally extremely intensive to run simulations for graphs of this size. The point here is that these simulations suggest that the GLRT is only able to detect at $\mu \approx m^{-1/6}$ and, therefore, does not achieve the optimal detection rate under the uniform prior on paths.

TABLE 1
*Value of the mean shift giving a detection rate of about 95% when using the Bayes test, the strip statistic test with width $B = 2\sqrt{m}$ and the GLRT (uniform prior on paths)—one can compute $\mu_{.95}$ for the strip statistic for large values of $m$ since it is given analytically*

| $m$ | 1025 | 2049 | 4097 | 8193 | 16385 |
|---|---|---|---|---|---|
| $\mu_{0.95}$ (Bayes) | 0.37 | 0.31 | 0.26 | | |
| $\mu_{0.95}$ (strip) | 0.84 | 0.69 | 0.59 | 0.51 | 0.42 |
| $\mu_{0.95}$ (GLRT) | 0.46 | 0.40 | 0.36 | 0.33 | |

6.2. *Minimax detection.* We focus here on the increasing path $p$, where $p_i = i$, $0 \le i \le m-1$, as we believe this path to be the most challenging for the GLRT. In this section, we compare the performance of the GLRT with the weighted average statistic test (WAS) defined in (2.12).

Recall that the WAS is distributed as $N(0, \lambda_m)$ under $H_0$ and as $N(\mu, \lambda_m)$ under $H_1$, regardless of the unknown path $[\lambda_m \sim (\log m)^{-1}]$. Thus, to achieve a power equal to 0.95 at the 5% significance level, we need $\mu \ge 2z_{0.95}\sqrt{\lambda_m}$, where $z_{0.95}$ is the 95% standard normal quantile. Some power curves for the WAS are graphed in Figure 5. We use simulations to graph similar curves

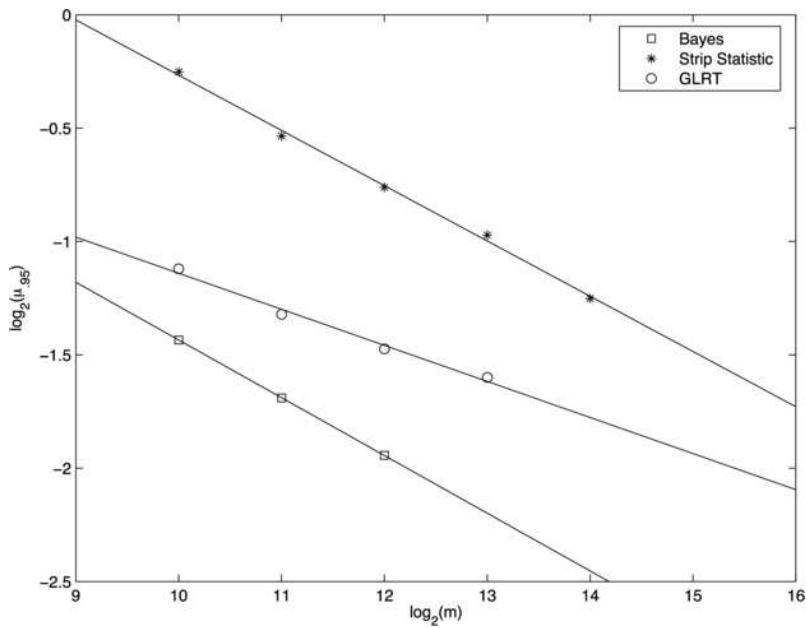

FIG. 3. *Comparison of the Bayes test, the strip statistic test and the GLRT under the uniform prior. The plot shows the value $\mu_{0.95}$ of the mean shift for which a given test achieves a 95% detection rate when the rate of false alarm is set at 5% as a function of the graph size $m$ (log-log scale).*

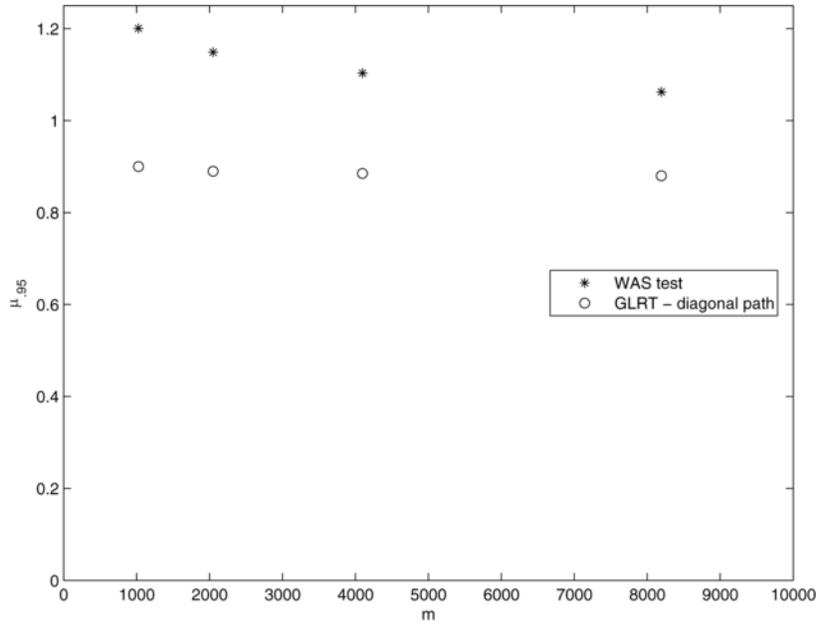

Fig. 4. *Comparison of the GLRT and the WAS when the anomalous path is the increasing path. The plot shows the value $\mu_{0.95}$ of the mean shift for which a given test achieves a 95% detection rate when the rate of false alarm is set at 5% as a function of the graph size $m$.*

Table 2
*Value of the mean shift giving a detection rate of about 95% when using the WAS test and the GLRT for detecting the increasing path—one can compute $\mu_{0.95}$ for the WAS for large values of $m$ since it is given analytically*

| $m$ | 1025 | 2049 | 4097 | 8193 | 16385 | 32769 |
|---|---|---|---|---|---|---|
| $\mu_{0.95}$ (WAS) | 1.20 | 1.15 | 1.10 | 1.06 | 1.03 | 0.99 |
| $\mu_{0.95}$ (GLRT) | 0.90 | 0.89 | 0.885 | 0.88 | | |

for the GLRT; see Figure 6. Each point is based on 1,000 realizations of the statistic.

While the power curves for the WAS tend to translate to the left, this does not seem to be the case for the GLRT. This might indicate that the detection threshold for the GLRT does not tend to zero as $m$ increases, just as in the case of the binary tree.

**7. Discussion.** Our paper leaves a number of open questions and invites several refinements. We briefly discuss some of these.

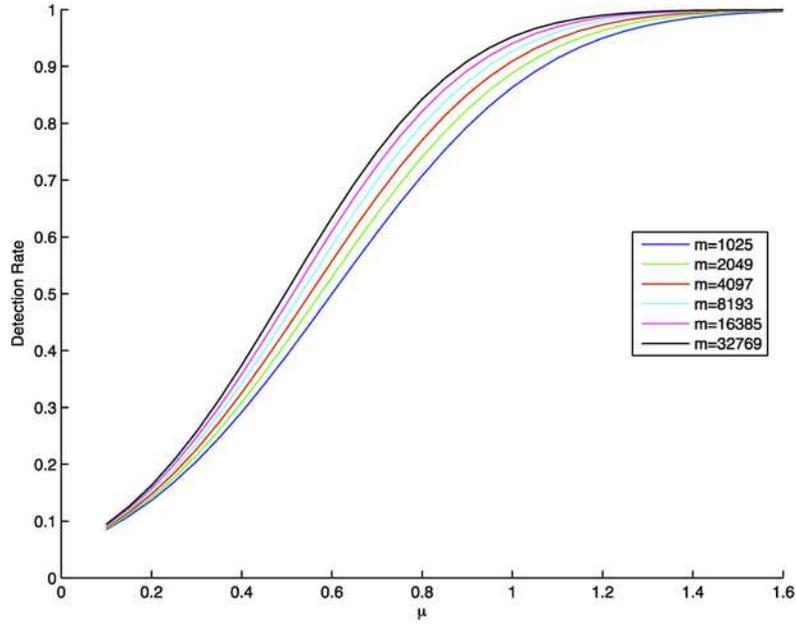

FIG. 5.  *Detection rate curves for the WAS statistic with $m = 1025$, 2049, 4097, 8193, 16385, 32769. As m increases, the curve moves to the left. The Type I error is set to 5%*

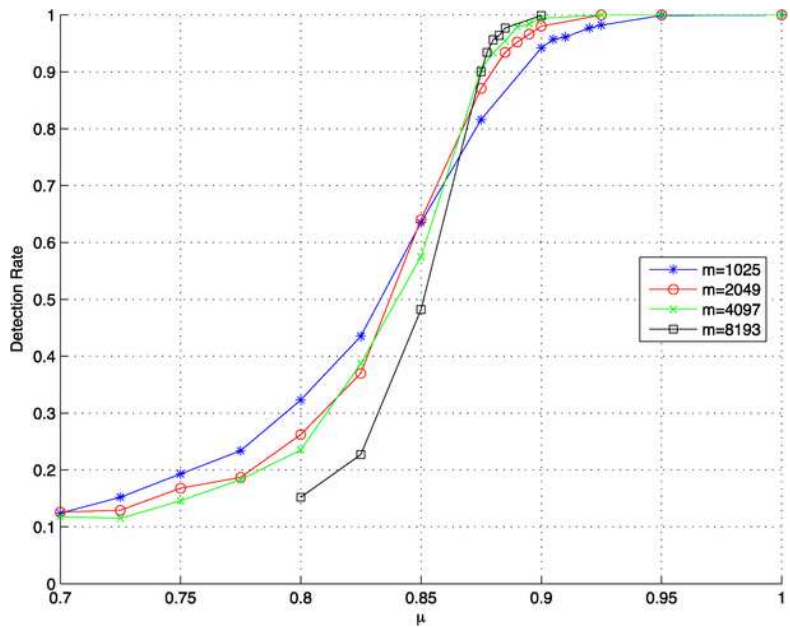

FIG. 6.  *Detection rate curves for the GLRT (increasing path). The probability of Type I error is set to 5%.*

7.1. *Sharpening the minimax detectability threshold in the two-dimensional regular lattice.* There is a gap between the upper and lower bounds in Theorem 1.1: the detection threshold for our estimator (2.12) is of order $\mu_m \sim (\log m)^{-1/2}$, but the priors we constructed showed nondetectability only when $\mu_m \sim (\log m)^{-1}$ (ignoring loglog factors). We do not see how to improve our prior to yield significantly better bounds and it seems that in any case, explicit priors of this family—as constructed in [17], for example— will not yield a lower bound obeying $\mu_m \gg (\log m)^{-3/4}$. It would be very interesting to understand this better and decide what is the actual rate of the detectability threshold.

With this in mind, we would like to emphasize that the test (2.12) used to prove the upper bound in Theorem 1.1 does not use the "continuity" of the path, only that it is known to be in the grid. That is, the test detects any sequence of the form $\{(i, p_i) : 0 \leq i \leq m - 1\}$ as long as $(i, p_i)$ is a vertex in the graph, provided, of course, that $\mu_m$ is of order $(\log m)^{-1/2}$. In fact, $(\log m)^{-1/2}$ turns out to be the minimax detection threshold when the set of vertices with positive mean is any sequence $(i, p_i)$ remaining in the grid. Indeed, the least favorable prior chooses the $(p_i)$ independently and uniformly at random in their respective range so that the number of crossings of two independent paths obeys

$$N_m = \sum_{1 \leq i \leq m} I_i,$$

where the $I_i$'s are independent with $\mathbf{P}(I_i = 1) = 1/i$ and $\mathbf{P}(I_i = 0) = 1 - 1/i$. The same argument as before shows that

$$\mathbf{E}_0(L_m - 1)^2 = \mathbf{E} e^{\mu_m^2 N_m} - 1 = \prod_{1 \leq i \leq m} \left(1 + \frac{e^{\mu_m^2} - 1}{i}\right) - 1,$$

which is easily shown to converge to zero when $\mu_m (\log m)^{1/2} \to 0$.

7.2. *Studying the GLRT on the two-dimensional regular lattice.* The GLRT may not be anywhere near optimal in the minimax sense. A indication of that can be deduced from work of Baik and Rains in [5], Section 4.4, and [6]. In the language of the current paper, they deal with the following problem: consider directed paths in the grid $\{(i, j) \in \mathbf{Z}^2 : j \leq i \leq m\}$. That is, starting from the origin $(0, 0)$, a path is a sequence of increments by 1 unit in the right or upward direction (this corresponds to a rotation of the regular graph considered in Theorem 1.1, with its lower half erased). Under $H_0$, all vertices are i.i.d. exponential random variables with parameter 1. Under $H_1$, the variables along the "diagonal path" [the path $(0, 0)$, $(1, 1)$, $(2, 2)$ and so on] are i.i.d. exponential with mean $1 + \mu$ (of course, in this situation, $H_1$ is asymptotically distinguishable from $H_0$ if $\mu > 0$, but this is of no concern in what follows). They consider the GLRT statistic $M_m$, which consists of the maximum partial sums among all possible directed paths connecting $(0,0)$ to $(m,m)$, and show that the limit distribution of (a properly rescaled version of) $M_m$ does not depend on $\mu$ as long as $\mu < 1$ (this follows from the geometric case treated in [5], Section 4.4). This hints that in that particular set-up, the GLRT is far from optimal since Section 2 shows that the minimax risk with respect to all possible directed path goes to zero for any $\mu > 0$. (Note that, strictly speaking, since the mode of convergence in [5] is weak convergence and not total variation, the results there hint, but do not imply, that the GLRT is not optimal.) Recently, Beffara and Sidoravicius (in a yet unpublished work) have analyzed the GLRT for the model considered in Theorem 1.1 (with exponential random variables), and their results seem to imply that the threshold for the GLRT is of order $o(1)$, in contrast with the case [5] treated by Baik and Rains.

Also of interest would be to study the power of the GLRT with a uniform prior on paths, where we suspect that the GLRT does not achieve the optimal threshold.

7.3. *Unknown starting location.* Throughout this paper, we assumed that under $H_1$, the unknown path starts at a known node (the origin). The same question can also be posed when the starting location is not known. For concreteness, consider the regular lattice as in Section 2 and allow the unknown path of length $m/2$ to start at any vertex in the collection $\{(i,j)\}_{i=0}^{m/2}$. Does there exist an asymptotically powerful test (in the minimax sense) for some sequence $\mu_m \to 0$? Similarly, we could also imagine having a square lattice $V_m = \{(i,j)\}$ with $0 \leq i \leq m-1$, $0 \leq j < 2m$ ($j$ has the parity of $i$ as before) and with edges $(i,j) \to (i+1, j+s)$, where $s = \pm 1$ and $j+s$ is understood modulo $2m$. If we know the starting location $(0,j)$ of the unknown path of length $m$, then this is the model problem discussed in Section 2. But studying this problem when we do not know the starting vertex is also of interest.

7.4. *Further refinements.* In this paper, we assumed that the node variables are independent and identically distributed and, clearly, one could address similar testing problems in far more general set-ups. Interesting extensions include situations in which the variables are correlated or in which the means along the unknown path are not all equal. Following up on the nonparametric signal detection problem, one could also imagine problems where the vector of means is not exactly sparse in the sense that it is zero away from the unknown path, but only rapidly decaying away from this path.

While this paper focuses on asymptotic properties of the detection problem, it is also of interest to develop test statistics with good finite sample size properties and we hope to report on our progress in a future publication.

SEARCHING FOR A TRAIL OF EVIDENCE IN A MAZE 317.5. *Other work.* While this paper was being written, N. Berger and Y. Peres described to us some of their own results, obtained independently, which address related problems and may answer some of the questions raised above.

## APPENDIX: PROOF OF LEMMA 2.3

To construct a stochastic process obeying (2.15), we follow [17] and let $S_n$ be the sum $S_n = \sum_{i=1}^n I_i$ with $\mathbf{P}(I_i = 1) = p_i$ and $\mathbf{P}(I_i = -1) = 1 - p_i$. Here, the $p_i$'s are stochastic (random environment) and defined by

$$p_i = 1/2 + p_i^{(1)} + p_i^{(2)} + \cdots,$$

where $(p_i^{(1)})$, $(p_i^{(2)})$, ..., are independent processes.

1. For each $i$ and $j$, the distribution of $p_i^{(j)}$ is uniform on $[-a_j, a_j]$.
2. The value $p_i^{(j)}$ is constant in $i$ for $i = 1, \ldots, 2^j$. At time $2^j + 1$, it switches to a new independent value, uniform on $[-a_j, a_j]$, which is kept until time $2 \times 2^j$, and so on.

Note that we need

(A.1) $$\sum_{j \geq 0} a_j < 1/2$$

for this to make sense so that the $p_i \in (0, 1)$. Finally, the $I_i$'s are independent, conditioned on the random environment $(p_i)$.

With this in place, Häggström and Mossel in [17], Proposition 3.1 showed that there exists a nearest-neighbor process $(S_n)$ obeying

(A.2) $$\text{PRE}_S(k) \leq \frac{C}{k a_{\lfloor \log_2(k/2) \rfloor}} \quad \text{for all } k = 1, 2, \ldots,$$

where $C = 4[C_1 + 1]$, with $C_1 = 2^{m_k} a_{m_k} \cdot \mathbf{P}(Y < \mathbf{E}Y/2)$, $m_k = \lfloor \log_2(k/2) \rfloor$ and $Y$ is a binomial random variable with $2^{m_k}$ trials and a probability of success equal to $a_{m_k}$. Since, for any binomial variable $Y_{n,p} \sim \text{Bin}(n, p)$,

$$np\mathbf{P}(Y_{n,p} < np/2) \leq \frac{4np\text{Var}(Y_{n,p})}{n^2 p^2} \leq 4,$$

$C_1 \leq 4$ and thus the constant $C \leq 20$.

As discussed earlier, this remark is of importance to us since we have used a sequence $(a_j)$ that depends explicitly on $m$.

**Acknowledgments.** Emmanuel Candès would like to thank the Centre Interfacultaire Bernoulli of the Ecole Polytechnique Fédérale de Lausanne for hospitality during June and July 2006. Emmanuel Candes and Ofer Zeitouni would like to thank Houman Owhadi for fruitful conversations. Ofer Zeitouni would like to thank Chris Hoffman and Yuval Peres for describing some of their unpublished results to him, and Jinho Baik for explaining to him the relevance of [5]; see the discussion in Section 7 and Vladas Sidoravicius for describing to him his work with Vincent Beffara.

## REFERENCES

[1] AHUJA, R. K., MAGNANTI, T. L. and ORLIN, J. B. (1993). *Network Flows. Theory, Algorithms and Applications.* Prentice Hall, Englewood Cliffs, NJ. MR1205775

[2] ARIAS-CASTRO, E., DONOHO, D. and HUO, X. (2005). Near-optimal detection of geometric objects by fast multiscale methods. *IEEE Trans. Inform. Theory* **51** 2402–2425. MR2246369

[3] ARIAS-CASTRO, E., DONOHO, D. and HUO, X. (2006). Adaptive multiscale detection of filamentary structures in a background of uniform random points. *Ann. Statist.* **34** 326–349. MR2275244

[4] BAHADUR, R. and RANGA RAO, R. (1960). On deviations of the sample mean. *Ann. Math. Statis.* **31** 1015–1027. MR0117775

[5] BAIK, J. and RAINS, E. M. (2001a). The asymptotics of monotone subsequences of involutions. *Duke Math. J.* **109** 205–281. MR1845180

[6] BAIK, J. and RAINS, E. M. (2001b). Symmetrized random permutations. In *Random Matrix Models and Their Applications* (P. Bleher and A. Its, eds.) 1–19. Cambridge Univ. Press. MR1842780

[7] BENJAMINI, I., PEMANTLE, R. and PERES, Y. (1998). Unpredictable paths and percolation. *Ann. Probab.* **26** 1198–1211. MR1634419

[8] BIGGINS, J. D. (1977). Martingale convergence in the branching random walk. *J. Appl. Probab.* **14** 25–37. MR0433619

[9] BUFFET, E., PATRICK, A. and PULÉ, J. V. (1993). Directed polymers on trees: A martingale approach. *J. Phys. A* **26** 1823–1834. MR1220795

[10] CANDES, E. J., CHARLTON, P. R. and HELGASON, H. (2006). Detecting highly oscillatory signals by chirplet path pursuit. Technical report, California Institute of Technology.

[11] COMETS, F., SHIGA, T. and YOSHIDA, N. (2004). Probabilistic analysis of directed polymers in random environments: A review. In *Stochastic Analysis on Large Scale Interacting Systems* 115–142. Math. Soc. Japan, Tokyo. MR2073332

[12] DERRIDA, B. and SPOHN, H. (1988). Polymers on disordered trees, spin glasses, and traveling waves. *J. Statist. Phys.* **51** 817–840. MR0971033

[13] DONOHO, D. L. and HUO, X. (2002). Beamlets and multiscale image analysis. In *Multiscale and Multiresolution Methods. Lecture Notes in Comput. Sci. Eng.* **20** 149–196. Springer, Berlin. MR1928566

[14] FELLER, W. (1957). The numbers of zeros and of changes of sign in a symmetric random walk. *Enseignement Math. (2)* **3** 229–235. MR0097864

[15] FELLER, W. (1968). *An Introduction to Probability Theory and Its Applications* I, 3rd ed. Wiley, New York. MR0228020

[16] GLAZ, J., NAUS, J. and WALLENSTEIN, S. (2001). *Scan Statistics*. Springer, New York. MR1869112

[17] HÄGGSTRÖM, O. and MOSSEL, E. (1998). Nearest-neighbor walks with low predictability profile and percolation in $2 + \epsilon$ dimensions. *Ann. Probab.* **26** 1212–1231. MR1640343
[18] HOFFMAN, C. (1998). Unpredictable nearest neighbor processes. *Ann. Probab.* **26** 1781–1787. MR1675075
[19] INGSTER, Y. I. and SUSLINA, I. A. (2003). *Nonparametric Goodness-of-fit Testing under Gaussian Models*. Springer, New York. MR1991446
[20] KULLDORFF, M. (1997). A spatial scan statistic. *Comm. Statist. Theory Methods* **26** 1481–1496. MR1456844
[21] LYONS, R. (1997). A simple path to Biggins' martingale convergence for branching random walk. In *Classical and Modern Branching Processes* (K. B. Athreya and P. Jagers, eds.) 217–221. Springer, New York. MR1601749
[22] PATIL, G. P., BALBUS, J., BIGING, G., JAJA, J., MYERS, W. L. and TAILLIE, C. (2004). Multiscale advanced raster map analysis system: Definition, design and development. *Environ. Ecol. Stat.* **11** 113–138. MR2086391
[23] PEMANTLE, R. (1995). Tree-indexed processes. *Statist. Sci.* **10** 200–213. MR1368099
[24] ZHONG, Y., JAIN, A. and DUBUISSON-JOLLY, M.-P. (2000). Object tracking using deformable templates. *IEEE Trans. Pattern Anal. Mach. Intell.* **2** 544–549.

E. ARIAS-CASTRO
DEPARTMENT OF MATHEMATICS
UNIVERSITY OF CALIFORNIA, SAN DIEGO
9500 GILMAN DRIVE
LA JOLLA, CALIFORNIA 92093-0112
USA
E-MAIL: eariasca@math.ucsd.edu

E. J. CANDÈS
H. HELGASON
APPLIED AND COMPUTATIONAL MATHEMATICS
CALIFORNIA INSTITUTE OF TECHNOLOGY
300 FIRESTONE, MAIL CODE 217-50
PASADENA, CALIFORNIA 91125
USA
E-MAIL: emmanuel@acm.caltech.edu
hannes@acm.caltech.edu

O. ZEITOUNI
DEPARTMENT OF MATHEMATICS
UNIVERSITY OF MINNESOTA
206 CHURCH STREET SE
MINNEAPOLIS, MINNESOTA 55455
USA
AND
DEPARTMENT OF MATHEMATICS
WEIZMANN INSTITUTE OF SCIENCE
POB 26
REHOVOT 76100
ISRAEL
E-MAIL: zeitouni@math.umn.edu